\documentclass[12pt]{amsart}

\usepackage{amssymb,amsfonts,amsmath,amsthm,epsfig,tikz,bm,color}
  
\def\Rm#1{\lowercase\expandafter{\romannumeral#1}}

   \def\Ga{\Gamma}

\def\a{\alpha} \def\b{\beta} \def\g{\gamma}

 \def\bbZ{{\Bbb Z}}  
   
  \def\nd{\mathrel{\bigm|\kern-.7em/}}

 \def\Aut{\hbox{\rm Aut}}

\def\Cay{\hbox{\rm Cay}}  
  
\def\qed{\hfill $\Box$} 
  \def\Arc{\hbox{\rm Arc}} 
 \def\bbZ{{\Bbb Z}} \def\bbQ{{\Bbb Q}} 
\def\Hol{\hbox{\rm Hol}} 

\newtheorem{thm}{Theorem}[section]

\newtheorem{lem}[thm]{Lemma}

\makeatletter \@addtoreset{equation}{section}

\newtheorem*{rmk1}{Remark}

\newtheorem{pro}[thm]{Proposition}

\newtheorem{prob}[thm]{Problem}
\def\pf{\noindent {\it Proof.\ }}
\def\qed{\ifmmode\square\else\nolinebreak\hfill
$\Box$\fi\par\vskip12pt}

\begin{document}

\title[ Normal Cayley digraphs of cyclic groups with CI-property]
{Normal Cayley digraphs of cyclic groups with CI-property}%

\author{Jin-Hua Xie, Yan-Quan Feng$^*$, Grigory Ryabov, Ying-Long Liu}%

\address{$^a$Department of Mathematics, Beijing Jiaotong University, Beijing, 100044, China\mbox{\hskip 3cm}  \linebreak
 \mbox{}\hskip 0.5cm $^b$Novosibirsk State University, Novosibirsk, Russia \mbox{\hskip 7cm}
  \linebreak
  \mbox{}\hskip 0.5cm $^c$Sobolev Institute of Mathematics, Novosibirsk, Russia}\mbox{\hskip 6cm}
  \email{ J.-H. Xie$^a$, jinhuaxie@bjtu.edu.cn; Y.-Q. Feng$^a$, yqfeng@bjtu.edu.cn; \linebreak
  \mbox{}\hskip 3.5cm G. Ryabov$^{b,c}$, gric2ryabov@gmail.com; Y.-L. Liu$^a$, 17121573@bjtu.edu.cn} %
\thanks{$^*$Corresponding author}


\begin{abstract}
A Cayley (di)graph $\Cay(G,S)$ of a group $G$ with respect to a subset $S$ of $G$ is called normal if the right regular representation of $G$ is a normal subgroup in the full automorphism group of $\Cay(G,S)$, and is called a CI-(di)graph if for every $T\subseteq G$, $\Cay(G,S)\cong \Cay(G,T)$ implies that there is $\sigma\in \Aut(G)$ such that $S^\sigma=T$. We call a group $G$ a NDCI-group if all normal Cayley digraphs of $G$ are CI-digraphs, and a NCI-group if all normal Cayley graphs of $G$ are CI-graphs, respectively. In this paper, we prove that a cyclic group of order $n$ is a NDCI-group if and only if $8\nmid n$, and is a NCI-group if and only if either $n=8$ or $8\nmid n$.
\end{abstract}

\maketitle
\qquad {\textsc k}{\scriptsize \textsc {eywords.}} {\footnotesize  Cayley digraph, DCI-group, CI-group, NDCI-group, NCI-group.}

\section{Introduction}


All (di)graphs considered in this paper are finite and simple, that is, no multiple arcs, edges or  loops, and all groups are finite. For a (di)graph $\Gamma$, denote by $V(\Gamma)$, $E(\Gamma)$, $\Arc(\Gamma)$ and $\Aut(\Gamma)$ the vertex set, edge set, arc set, and full automorphism group of $\Gamma$, respectively.

Let $G$ be a finite group and let $S$ be a subset of $G$ with $1\not\in S$. The {\em Cayley digraph $\Cay(G,S)$} of $G$ with respect to $S$ is defined to be the digraph with vertex set $G$ and arc set $\{(g,sg)\ |\ g\in G, s\in S\}$. A Cayley digraph $\Cay(G,S)$ is called a {\em Cayley graph} of $G$ with respect to $S$ if $S$ is an inverse-closed set, that is, $S=S^{-1}:=\{x^{-1}\mid x\in S\}$, and in this case, $\Cay(G,S)$ is viewed as an undirected graph by identifying two arcs $(u,v)$ and $(v,u)$ as an edge $\{u,v\}$. Thus, Cayley graph is a special case of Cayley digraph.

For a given $g\in G$, the right multiplication $R(g)$ is a permutation on $G$ such that $x^{R(g)}=xg$ for all $x\in G$. Then the permutation group $R(G)=\{R(g)\ |\ g\in G\}$ on $G$ is called the {\em right regular representation} of $G$, and it is a regular subgroup of $\Aut(\Cay(G,S))$, so that $\Cay(G,S)$ is vertex-transitive. In fact, Sabidussi~\cite{Sabidussi} proved that a digraph $\Gamma$ is a Cayley digraph of a group $G$ if and only if $\Aut(\Gamma)$ contains a regular subgroup isomorphic to $G$. Furthermore, a Cayley digraph $\Cay(G,S)$ has out-valency $|S|$, and is connected if and only if $\langle S\rangle=G$, namely, $S$ generates $G$.

Let $G$ be a group. Two Cayley (di)graphs $\Cay(G,S)$ and $\Cay(G,T)$ of $G$ are called Cayley isomorphic if $S^{\sigma}=T$ for some $\sigma\in \Aut(G)$. It is easy to see that Cayley isomorphic Cayley (di)graphs are isomorphic as (di)graphs. However, the converse is not true, and there are isomorphic Cayley (di)graphs which are not Cayley isomorphic.
A subset $S$ of $G$ with $1\not\in S$ is called a {\em CI-subset} if for $T\subseteq G$, $\Cay(G,S)\cong \Cay(G,T)$ implies that they are Cayley isomorphic, that is, $T=S^{\sigma}$ for some $\sigma\in \Aut(G)$, and in this case, $\Cay(G,S)$ is said to be a {\em CI-digraph} or a {\em CI-graph} when it is a digraph or graph respectively, where CI represents Cayley isomorphic. A group $G$ is said to be a {\em DCI-group} or a {\em CI-group} if all Cayley digraphs or graphs of $G$ are CI-digraphs or CI-graphs, respectively.

A Cayley (di)graph $\Cay(G,S)$ is called {\em normal} if the right regular representation $R(G)$ of $G$ is a normal subgroup of $\Aut(\Cay(G,S))$. We call a group $G$ a {\em NDCI-group} or a {\em NCI-group} if all normal Cayley digraphs or graphs of $G$ are CI-digraphs or CI-graphs, respectively. Obviously, a DCI-group is a NDCI-group and a CI-group is a NCI-group. In this paper, we determine  cyclic NDCI-groups and cyclic NCI-groups.

\begin{thm}\rm\label{mainth}
A cyclic group of order $n$ is a $\rm{NDCI}$-group if and only if $8\nmid n$. A cyclic group of order $n$ is a $\rm{NCI}$-group if and only if either $n=8$ or $8\nmid n$.
\end{thm}

A lot of efforts have been made to achieve classifications of DCI-groups and CI-groups. We do not intend to give a full account on their study in general, and interested reader is referred to \cite{Kov,KR,C.H.Li,LLP,MS}, and in particular to \cite{AN,CL,D2,Feng,HM,J.M,Mu,N,So,Sp,Sp1} for elementary abelian groups, where many authors like Li, Kov\'{a}cs, Spiga et al. gave a lot of contributions.
Here we only outline classifications of cyclic DCI-groups and CI-groups, which originated from a conjecture proposed by \'Ad\'am~\cite{Adam}: every finite cyclic group is a CI-group. Though the conjecture was disproved by Elspas and Turner~\cite{Elspas}, it started the study of classifications of cyclic DCI-groups and CI-groups, which had lasted for 30 years. It was known that the \'Ad\'am's conjecture is true for the cyclic group of order $n$ in each of the following cases ($p$ and $q$ are distinct primes): $n=p$~\cite{Djokovic,Elspas,Turner}; $n=2p$~\cite{Babai}; $n=pq$~\cite{Alspach,Godsil}; $n=4p$ with $p>2$~\cite{Godsil}; $(n, \phi(n))=1$~\cite{Palfy} where $\phi$ is the Euler's function. Finally, Muzychuk~\cite{Muzychuk,M.Muzychuk} put the last piece into the puzzle and proved that a cyclic group of order $n$ is a DCI-group if and only if $n=mk$ where $m=1,2,4$ and $k$ is odd square-free, and is a CI-group if and only if either $n=8,9,18$, or $n=mk$ where $m=1,2,4$ and $k$ is odd square-free. This is one of the remarkable achievements on the study of DCI-groups and CI-groups.

Normality of Cayley (di)graphs has also been investigated intensively. Godsil~\cite{Godsil1} proved that if $\Cay(G,S)$ is normal then $\Aut(\Cay(G,S))$ is a semidirect product of $R(G)$ by the subgroup $\Aut(G,S):=\{\a\in\Aut(G)\ |\ S^\a=S\}$ (also see Proposition~\ref{N_AUT}), and hence $\Aut(\Cay(G,S))$ is completely determined by $\Aut(G)$. Wang et al.~\cite{C.Q.Wang} obtained that every finite group $G$ has a normal Cayley graph unless $G\cong C_4 \times C_2$ or $G \cong \bbQ_8 \times C_2^m$, where $C_n$ stands a cyclic group of order $n$ and $\bbQ_8$ is the Quaternion group. The normality of Cayley graphs of cyclic groups of order a prime and of groups of order twice a prime was solved by Alspach~\cite{Alspach1} and Du et al.~\cite{DuSF}, respectively. Dobson~\cite{D} determined all non-normal Cayley graphs of order a product of two distinct primes, and Dobson and Witte~\cite{DD} determined all non-normal Cayley graphs of order a prime-square. For normality of Cayley graphs of finite simple groups, we refer the reader to \cite{Li,FPW,FLWX,FLX,FLuWX,XFWX1,XFWX2,FMW,DFZ,DF,PYL,YFZC},
and for  some results on normality of Cayley digraphs, one may see \cite{Baik,Feng1,Feng2,Zhou}.
Based on these results, Xu~\cite{Xu} conjectured that almost all connected Cayley graphs are normal.

Li~\cite{C.H.Li} constructed some normal Cayley digraphs of cyclic groups of $2$-power order which are not CI-graphs, and proposed the following problem: Characterize normal Cayley digraphs which are not CI-graphs. Motivated in part by this and classifications of cyclic DCI-groups and CI-groups, we would like to propose the following problem.

\begin{prob}\label{prob}
Classify finite NDCI-groups and NCI-groups.
\end{prob}

Theorem~\ref{mainth} completely solves the problem for cyclic groups. It seems that Problem~\ref{prob} is easier than the classification of finite DCI-groups and CI-groups, but still quite difficult. However, solving Problem~\ref{prob} should be helpful for the classification of finite DCI-groups and CI-groups, which is very difficult.

All the notations and terminologies used throughout this paper are
standard. For group and graph theoretic concepts which are not
defined here we refer the reader to \cite{Biggs,Dixon,Rotman,Wi}.

\section{Preliminaries}
In this section, we give some basic concepts and facts that will be needed later.
Let $\Cay(G,S)$ be a Cayley digraph of a group $G$ with respect to $S$. Recall that  $\Aut(G,S)=\{ \alpha\in \Aut(G)~|~S^{\alpha}=S\}$. Then $\Aut(G,S)$ is a subgroup of $\Aut(\Cay(G,S))$, and in fact, it is a subgroup of the stabilizer $(\Aut(\Cay(G,S)))_1$ of $1$ in $\Aut(\Cay(G,S))$. Godsil~\cite{Godsil1} proved that the normalizer of $R(G)$ in $\Aut(\Cay(G,S))$ is a semiproduct of $R(G)$ by $\Aut(G,S)$, and then one may prove the following proposition easily (see \cite[Propositions 1.3 and 1.5]{Xu}).

\begin{pro}\rm\label{N_AUT} Let $\Cay(G,S)$ be a Cayley digraph of a group $G$ with respect to $S$ and let $A=\Aut(\Cay(G,S))$. Then $N_A(R(G))=R(G)\rtimes\Aut(G,S)$ and
$\Cay(G,S)$ is normal if and only if $A_1=\Aut(G,S)$.
\end{pro}

Babai~\cite{Babai} obtained the well-known criterion for a Cayley digraph to be a CI-digraph: a Cayley digraph $\Cay(G,S)$ is a CI-digraph if and only if all regular subgroups of $\Aut(\Cay(G,S))$ isomorphic to $G$ are conjugate. This was also proved by Alspach and Parsons~\cite{Alspach}. Based on this criterion, the follow proposition is straightforward (also see \cite[Corollary 6.9]{C.H.Li}).

\begin{pro}\rm\label{CI-graph-prop}
Let $\Cay(G,S)$ be a normal Cayley digraph of a group $G$ with respect to $S$. Then $\Cay(G,S)$ is a $\rm{CI}$-digraph if and only if $R(G)$ is the unique regular subgroup of $\Aut(\Cay(G,S))$ which is isomorphic to $G$.
\end{pro}
Let $G$ be a group and let $g\in G$. We denote by $o(g)$ the order of $g$ in $G$. Li~\cite[Example 6.10]{C.H.Li} constructed some normal Cayley digraphs of cyclic groups of $2$-power order which are not CI-digraphs.

\begin{pro}\rm \label{pro-no-ci}
Let $G=\langle a\rangle$ with $o(a)=2^r\geq 8$ and let $S=\{a,a^2,a^{2^{r-1}+1}\}$. Then $\Cay(G,S)$ is a normal Cayley digraph but not a $\rm{CI}$-digraph.
\end{pro}

The classifications of cyclic DCI-groups and CI-groups were finally completed by Muzychuck \cite{Muzychuk,M.Muzychuk} (also see \cite[Theorem 7.1]{C.H.Li}).

\begin{pro}\rm\label{AdNe-DCI-group}
 A cyclic group of order $n$ is a DCI-group if and only if $n =k, 2k$ or $4k$ where $k$ is odd square-free, and a cyclic group of order $n$ is a CI-group if and only if either $n\in \{8,9,18\}$ or $n=k, 2k$ or $4k$ where $k$ is odd  square-free.
\end{pro}

Let $G$ be a group and let $G_i\leq G$ for $1\leq i\leq n$, with
$$G=G_1\times G_2\times\cdots\times G_n,$$ that is, $G$ is a direct product of subgroups $G_1,G_2,\cdots,G_n$. Then $\Aut(G_i)$ can be naturally extended to a group of automorphisms of $G$: for every $\a_i\in \Aut(G_i)$, define
\begin{equation}\label{eq1}
(g_1\cdots g_{i-1} g_ig_{i+1}\cdots g_n)^{\a_i}=g_1\cdots g_{i-1}g_i^{\a_i}g_{i+1}\cdots g_n,\mbox{where } g_j\in G_j,\ 1\leq j\leq n.
\end{equation}
For not making the notation too cumbersome, we still denote  by $\Aut(G_i)$ this extended group of automorphisms of $G$, which will be used throughout the paper. In particular, we have the following proposition for cyclic groups (see \cite[Theorem 7.3]{Rotman}).

\begin{pro} \rm\label{AGCG}  For a positive integer $n$, let $n=\prod_{i=1}^m{p_i^{r_i}}$ be the distinct prime factorization of $n$, and let $C_n=C_{p_1^{r_1}}\times C_{p_2^{r_2}}\times\cdots\times C_{p_m^{r_m}}$. Then $\Aut(C_n)=\Aut(C_{p_1^{r_1}})\times \Aut(C_{p_2^{r_2}})\times\cdots\times\Aut(C_{p_m^{r_m}})$. Furthermore, $\Aut(C_2)=1$, $\Aut(C_4)\cong C_2$, $\Aut(C_{2^n})\cong C_2\times C_{2^{n-2}}$ for $n\geq 2$, and $\Aut(C_{p^n})\cong C_{(p-1)p^{n-1}}$ for an odd prime $p$.
\end{pro}

Let $\Ga$ be a digraph and let $T,L\subseteq V(\Ga)$ with $T\cap L=\emptyset$. Denote by $[T]$ the induced sub-digraph of $T$ in $\Gamma$, that is, the digraph with the vertex set $T$ and with an arc $(u,v)\in \Arc([T])$, $u,v\in T$, whenever $(u,v)\in \Arc(\Gamma)$, and denote by $[T,L]$ the induced bipartite sub-digraph with vertex set $T\cup L$ and all arcs from $T$ to $L$ in $\Gamma$.

The following proposition describes a family of non-normal Cayley digraphs, and its main idea comes from the so-called generalized wreath circulants~\cite{BDM} and generalized wreath products of Cayley schemes~\cite[Section~3.4.3]{CP}.

\begin{pro} \rm\label{General-normal}  Let $G$ be a finite group, and let $1\not\in S\subset G$, $1<H\leq K<G$ and $H\unlhd G$. Assume that $S\backslash K$ is a union of some cosets of $H$ in $G$ and that there exist $x\not\in K$ and $y\in H$ such that $y^x\not=y^{-1}$. Then the Cayley digraph  $\Cay(G,S)$ is non-normal.
\end{pro}

\pf Write $\Ga=\Cay(G,S)$ and $A=\Aut(\Ga)$. Note that $R(G)\leq A$.

Take $u,v\in G$ with $Ku\not=Kv$. Note that for $g\in G$, the out-neighbourhood of $g$ in $\Ga$ is $\{sg\ |\ s\in S\}$. Since $S\backslash K$ is a union of some cosets of $H$ in $G$, for all cosets $Hk_1u$ of $H$ in $Ku$ and all cosets $Hk_2v$ in $Kv$, where $k_1,k_2\in K$, the normality of $H$ in $G$ implies that the induced bipartite sub-digraph $[Hk_1u, Hk_2v]$ is isomorphic to either the empty graph with $2|H|$ isolated vertices, or the complete bipartite digraph $\vec{K}_{|H|,|H|}$ of order $2|H|$.
By the arbitrariness of $u$ and $v$, if an automorphism of the induced sub-digraph $[Ku]$ of $Ku$ in $\Ga$ fixes every coset of $H$ in $Ku$ setwise, then it can be extended to an automorphism of $\Gamma$ fixing $G\backslash Ku$ pointwise.

For every $h\in H$, define $\bar{h}$ to be the permutation on $V(\Ga)$ such that $\bar{h}$ fixes $G\backslash Kx$ pointwise and $t^{\bar{h}}=t^{R(h)}=th$ for all $t\in Kx$, where $x\not\in K$ by hypothesis. Then $\bar{H}:=\{\bar{h}\ |\ h\in H\}$ is a permutation group on $V(\Ga)$. Clearly, the restriction of $\bar{H}$ on $Kx$ is the same as $R(H)$, and since $R(H)\leq A$ and $H\unlhd G$,  this restriction is a group of automorphisms of the induced sub-digraph $[Kx]$ of $Kx$ in $\Ga$ fixing every coset of $H$ setwise, and then the above paragraph implies that $\bar{H}\leq A$. Since $Kx\not=K$, we have $\bar{H}\leq A_1$, where $A_1$ is the stabilizer of $1$ in $A$.

Suppose to the contrary that $\Ga$ is normal. By Proposition~\ref{N_AUT}, $A_1=\Aut(G,S)$ and hence $\bar{H}\leq \Aut(G)$. Since $Kx^2\not=Kx$, $\bar{H}$ fixes $Kx^2$ pointwise. Since $y\in H$ by hypothesis, we have $\bar{y}\in \bar{H}$, and hence $(x^2)^{\bar{y}}=x^2$ and $x^{\bar{y}}=xy$ as $x\in Kx$. On the other hand, since $\bar{H}\leq \Aut(G)$ we obtain  $(x^2)^{\bar{y}}=x^{\bar{y}}x^{\bar{y}}=xyxy$.
It follows that $xyxy=x^2$, that is, $y^x=y^{-1}$, contrary to the hypothesis $y^x\not=y^{-1}$. Thus, $\Ga$ is non-normal.\qed

\begin{rmk1}\rm
The fact that $\bar{H}\leq A_1$ from the proof of Proposition~\ref{General-normal} was proved in~\cite[Theorem~3.4.21]{CP}  in the language of Cayley schemes.
\end{rmk1}

\section{Proof of Theorem~\ref{mainth}}

We first prove a crucial lemma for the proof of Theorem~\ref{mainth}.

\begin{lem}  \rm\label{non-normal}
For a positive integer $n$, let $n=\prod_{i=1}^m{p_i^{r_i}}$ be the distinct prime factorization of $n$, and let $\Cay(C_n,S)$ be a Cayley digraph of the cyclic group $C_n=C_{p_1^{r_1}}\times C_{p_2^{r_2}}\times\cdots\times C_{p_m^{r_m}}$. Assume that $p_t$ is an odd prime and that $\Aut(C_n,S)$ contains an element of order $p_t$ in $\Aut(C_{p_t^{r_t}})$ for some $1\leq t\leq m$. Then $\Cay(C_n,S)$ is non-normal.
\end{lem}

\pf Let  $\Ga=\Cay(C_n,S)$ and $A=\Aut(\Ga)$. Write $C_n=C_{p_1^{r_1}}\times C_{p_2^{r_2}}\times\cdots\times C_{p_m^{r_m}}=\langle a_1\rangle\times\langle a_2\rangle\times\dots\times\langle a_m\rangle$ with $o(a_i)=p_i^{r_i}$.  Then $C_{p_i^{r_i}}=\langle a_i\rangle$ for each $1\leq i\leq m$, and $C_n=\langle a_1a_2\ldots a_m\rangle$. Furthermore, each $x\in C_n$ can be uniquely factorized as $x=a_1^{s_1}a_2^{s_2}\cdots a_i^{s_i}\cdots a_m^{s_m}$, and we call $a_i^{s_i}$ the $i$-th component of $x$, denoted by $x_i$, that is, $x_i=a_i^{s_i}\in \langle a_i\rangle=C_{p_i^{r_i}}$.

Since $p_t$ is odd, $\Aut(C_{p_t^{r_t}})\cong C_{(p_t-1)p_t^{r_t-1}}$, and since $\Aut(C_n,S)$ contains an element of order $p_t$ in $\Aut(C_{p_t^{r_t}})$, we have $r_t\geq 2$. Note that the map  $a_t\mapsto a_t^{p_r^{r_t-1}+1}$ induces an automorphism of order $p_t$ of $C_{p_t^{r_t}}$.
By Eq~(\ref{eq1}), an element $\delta$ of order $p_t$ in $\Aut(C_{p_t^{r_t}})$ can be defined by
$$(a_1^{j_1}a_2^{j_2}\dots a_{t-1}^{j_{t-1}}a_t^{j_t}a_{t+1}^{j_{t+1}}\dots a_m^{j_m})^\delta= a_1^{j_1}a_2^{j_2}\dots a_{t-1}^{j_{t-1}}(a_t^{j_t})^{p_t^{r_t-1} + 1}a_{t+1}^{j_{t+1}}\dots a_m^{j_m}$$ for every $a_k^{j_k}\in C_{p_k^{r_k}}$, or alternatively by $$x^\delta=x(x_t)^{p_t^{r_t-1}} \mbox{ for every } x\in C_n.$$ The uniqueness of subgroup of order $p_t$ in $\Aut(C_{p_t^{r_t}})$ implies that $\langle \delta\rangle\leq \Aut(C_n,S)$. Let
$$K=\langle a_1,a_2,\cdots,a_{t-1},a_t^{p_t},a_{t+1},\cdots,a_m\rangle, \mbox{ and } H=\langle a_t^{p_t^{r_t-1}}\rangle.$$

Then $|C_n:K|=p_t$ and $|H|=p_t$. Since $r_t\geq 2$, $K$ has a subgroup of order $p_t$, and since $C_n$ has a unique subgroup of order $p_t$, we have $1<H\leq K<C_n$. By taking $x=a_t$ and $y=a_t^{p_t^{r_t-1}}$, we have $y^x\not=y^{-1}$. To finish the proof, by Proposition~\ref{General-normal}, it suffices to show that $S\backslash K$ is a union of some cosets of $H$ in $G$.

Take $y\in S\backslash K$. Then $y_t$ has order $p_t^{r_t}$ as $C_n=\prod_{i=0}^{p_t-1}a_t^iK$, implying $y_t^{p_t^{r_t-1}}\not=1$. Since $y^\delta=y(y_t)^{p_t^{r_t-1}}$, $\langle \delta \rangle$ is transitive on $yH$, and since $\langle \delta \rangle\leq \Aut(C_n,S)\leq A$, we have $yH\subseteq S\backslash K$. By the arbitrariness of $y$,  $S\backslash K$ is a union of some cosets of $H$ in $G$, as required. Thus, $\Gamma$ is non-normal.\qed

\medskip
Now we construct normal Cayley graphs on cyclic groups which are not CI-graphs.

\begin{lem}\rm\label{thm-no-ci}
Let $m$ and $s$ be positive integers such that $s\geq 3$, $(2,m)=1$ and ${2^s}m\neq 8$. Let $C_{2^sm}=\langle a \rangle\times\langle b\rangle$ with $o(a)=2^s$ and $o(b)=m$, and $S=\{(ab)^{\pm1},(a^2b^2)^{\pm1},(a^{2^{s-1}+1}b)^{\pm1}\}$. Then $\Cay(C_{2^sm},S)$ is a normal Cayley graph but not a CI-graph.
\end{lem}

\pf Let $\Gamma=\Cay(C_{2^sm},S)$ and $A=\Aut(\Gamma)$. Since $S=S^{-1}$, $\Gamma$ is a Cayley graph. For $u,v\in V(\Gamma)$, denote by $d(u,v)$ the distance between $u$ and $v$ in $\Gamma$. Write $\Gamma_k(u)=\{w \ |\ d(u,w)=k\}$, the $k$-step neighborhood of the vertex $u$ in $\Gamma$.

Let $A_1$ be the stabilizer $1$ in $A$, and let  $A^*_1=\{\a\in A_1\ |\ u^\a=u \mbox{ for every } u\in \Gamma_1(1)=S\}$, namely, the subgroup of $A_1$ fixing $S$ pointwise. Define $\a$ to be the automorphism of $C_{2^sm}$ induced by $a\mapsto a^{-1}$ and $b\mapsto b^{-1}$, and $\b$ to be the automorphism of $C_{2^sm}$ induced by $a\mapsto a^{2^{s-1}+1}$ and $b\mapsto b$. It is easy to check that $C_2\times C_2\cong\langle \a,\b \rangle\leq\Aut(G,S)\leq A_1$.

\begin{figure}[htb]
\begin{center}
\begin{picture}(330,145)(20,30)
\thicklines

\put(180,160){\circle*{5}}
\put(20,100){\circle*{5}}
\put(100,100){\circle*{5}}
\put(260,100){\circle*{5}}
\put(340,100){\circle*{5}}
\put(60,80){\circle*{5}}
\put(300,80){\circle*{5}}
\put(20,40){\circle*{5}}
\put(60,40){\circle*{5}}
\put(100,40){\circle*{5}}
\put(140,40){\circle*{5}}
\put(180,40){\circle*{5}}
\put(220,40){\circle*{5}}
\put(260,40){\circle*{5}}
\put(300,40){\circle*{5}}
\put(340,40){\circle*{5}}

\qbezier(180,160)(20,140)(20,100)
\qbezier(180,160)(110,140)(100,100)
\qbezier(180,160)(250,140)(260,100)
\qbezier(180,160)(320,140)(340,100)
\qbezier(180,160)(60,140)(60,80)
\qbezier(180,160)(300,140)(300,80)

\put(20,100){\line(1,0){80}}
\put(20,100){\line(2,-1){40}}
\put(20,100){\line(0,-1){60}}
\put(20,100){\line(2,-3){40}}
\put(20,100){\line(4,-3){80}}
\put(60,80){\line(0,-1){40}}
\put(60,80){\line(2,-1){80}}
\put(60,80){\line(3,-1){120}}
\put(60,80){\line(1,0){240}}
\put(100,100){\line(1,0){160}}
\put(100,100){\line(0,-1){60}}
\put(100,100){\line(2,-3){40}}
\put(100,100){\line(4,-3){80}}
\put(260,100){\line(1,0){80}}
\put(260,100){\line(0,-1){60}}
\put(260,100){\line(-2,-3){40}}
\put(260,100){\line(-4,-3){80}}
\put(300,80){\line(0,-1){40}}
\put(300,80){\line(-2,-1){80}}
\put(300,80){\line(-3,-1){120}}
\put(340,100){\line(0,-1){60}}
\put(340,100){\line(-2,-1){40}}
\put(340,100){\line(-2,-1){120}}
\put(340,100){\line(-2,-3){40}}
\put(340,100){\line(-4,-3){80}}

\put(178,165){${\small 1}$}
\put(24,102){${\small a^2b^2}$}
\put(104,102){${\small a^{2^{s-1}+1}b}$}
\put(204,102){${\small a^{2^{s-1}-1}b^{-1}}$}
\put(297,102){${\small a^{-2}b^{-2}}$}
\put(63,82){${\small ab}$}
\put(265,82){${\small a^{-1}b^{-1}}$}
\put(10,26){${\small a^4b^4}$}
\put(43,26){${\small a^3b^3}$}
\put(68,26){${\small a^{2^{s-1}+3}b^3}$}
\put(115,26){${\small a^{2^{s-1}+2}b^2}$}
\put(165,26){${\small a^{2^{s-1}}}$}
\put(188,26){${\small a^{2^{s-1}-2}b^{-2}}$}
\put(242,26){${\small a^{2^{s-1}-3}b^{-3}}$}
\put(295,26){${\small a^{-3}b^{-3}}$}
\put(335,26){${\small a^{-4}b^{-4}}$}
\end{picture}
\end{center}
\caption{The induced subgroup $[\Gamma_0(1)\cup \Gamma_1(1)\cup \Gamma_2(1)]$ in $\Gamma$.}
\label{figure2}
\end{figure}
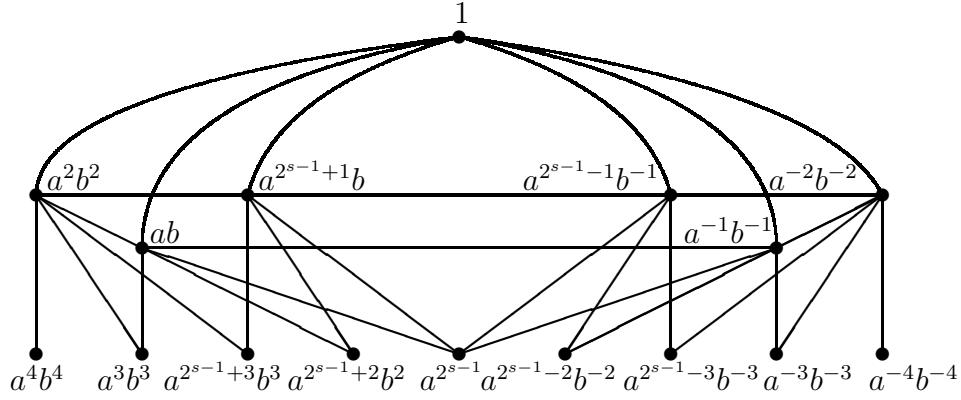

Clearly, $\Gamma_2(1)=\{(a^4b^4)^{\pm1},(a^3b^3)^{\pm1},(a^{2^{s-1}+3}b^3)^{\pm1},(a^{2^{s-1}+2}b^2)^{\pm1},a^{2^{s-1}}\}$.
Since $s\geq 3$, $(2,m)=1$ and ${2^s}m\neq 8$, we derive that $|\Gamma_2(1)|=9$, and the induced subgroup $[\Gamma_0(1)\cup \Gamma_1(1)\cup \Gamma_2(1)]$ can be depicted as Figure~\ref{figure2}.

Note that $A_1^*$ fixes $1$ and $S$ pointwise. For every two distinct vertices $u,v\in \Gamma_2(1)$, from Figure~\ref{figure2} we have that $\Gamma_1(u)\cap S\not=\Gamma_1(v)\cap S$, that is, every two distinct vertices in $\Gamma_2(1)$ cannot have the same neighbours in $\Gamma_1(1)$. Since $A_1^*$ fixes $\Gamma_0(1)\cup \Gamma_1(1)$ pointwise and $\Gamma_2(1)$ setwise, $A_1^*$ fixes $u$ and $v$ and hence $A_1^*$ fixes $\Gamma_2(1)$ pointwise. By the transitivity of $A$ on $V(\Gamma)$, $A_w^*$ fixes $\Gamma_2(w)$ pointwise, for every $w\in V(\Gamma)$. Since $\langle S\rangle=C_{2^sm}$, $\Gamma$ is connected, and an easy inductive argument on $k$ gives rise to that $A_1^*$ fixes $\Gamma_k(1)$ pointwise for every positive integer $k$. Thus, $A_1^*$ fixes each vertex of $\Gamma$, that is, $A_1^*=1$, and in particular, $A_1$ acts faithfully on $S$.

Again by Figure~\ref{figure2}, the induced subgraph $[S]$ is a cycle of length $6$ and hence $\Aut([S])\cong D_{12}$, the dihedral group of order $12$. Since $A_1$ is faithful on $S$, we may assume that $A_1\leq \Aut([S])$, and since $C_2\times C_2\cong\langle \a,\b \rangle\leq\Aut(G,S)\leq A_1$, we have either $A_1=\Aut([S])\cong D_{12}$, or $A_1=\Aut(G,S)$.

Suppose that $A_1=\Aut([S])\cong D_{12}$. Then $A_1$ has a unique Sylow $3$-subgroup of order $3$, and since $[S]$ is a cycle of length $6$, the unique Sylow $3$-subgroup can be generated by an element of order $3$, say $\g$, such that $(ab)^\g=a^{2^{s-1}+1}b$, $(a^{2^{s-1}+1}b)^\g=a^{-2}b^{-2}$ and $(a^{-2}b^{-2})^\g=ab$. From Figure~\ref{figure2}, $a^{2^{s-1}}$ is the unique vertex in $\Gamma_2(1)$ that has four neighbours in $\Gamma_1(1)$, which yields that $A_1$ fixes $a^{2^{s-1}}$, namely,  $(a^{2^{s-1}})^\g=a^{2^{s-1}}$. Since $\{ab,a^{2^{s-1}}\}\in E(\Gamma)$, we have $\{ab,a^{2^{s-1}}\}^{\g^{-1}}\in E(\Gamma)$, that is, $\{a^{-2}b^{-2},a^{2^{s-1}}\}\in E(\Gamma)$, which is not true by  Figure~\ref{figure2}.

It follows that $A_1=\Aut(G,S)$. By Proposition~\ref{N_AUT}, $\Gamma$ is normal. To finish the proof, we are left to show that $\Gamma$ is not a $\rm{CI}$-graph.

Let $H=\langle R(ab)\b\rangle$. Since $\b^{-1}R(ab)\b=R((ab)^\b)=R(a^{2^{s-1}+1}b)$ and $\b^2=1$, we have $$(R(ab)\b)^2=R(ab)R(ab)^\b=R(ab)R(a^{2^{s-1}+1}b)=R(a^{2^{s-1}+2}b^2),$$ and since $o(a^{2^{s-1}+2}b^2)=2^{s-1}m$, we obtain that $o(R(ab)\b)$ is even and $o(R(ab)\b)=(2,o(R(ab)\b))o((R(ab)\b)^2)=2^sm$, that is, $H\cong C_{2^sm}$.
Note that $\langle R(a^{2^{s-1}+2}b^2)\rangle$ is semiregular on $V(\Gamma)$ with two orbits of length $2^{s-1}m$, that is, $\langle a^{2^{s-1}+2}b^2\rangle=\langle a^2b\rangle$ and $a\langle a^2b\rangle$. Since $1^{R(ab)\b}=(ab)^\b=a^{2^{s-1}+1}b\in a\langle a^2b\rangle$ and $\langle R(a^{2^{s-1}+2}b^2)\rangle$ is a normal subgroup of index $2$ in $H$, $R(ab)\b$ interchanges the two orbits $\langle a^2b\rangle$ and $a\langle a^2b\rangle$, implying that $H$ is transitive on $V(\Gamma)$. Since $|H|=|V(\Gamma)|$, $H$ is regular on $V(\Gamma)$, and  since $R(C_{2^{s}m})\not=H$, $\Gamma$ is not a {\rm CI}-graph by Proposition~\ref{CI-graph-prop}. \qed

\medskip
Let $G$ be a finite group and let $p$ be a prime. Denote by $G_p$ a Sylow $p$-subgroup of $G$. The right multiplication $R(G)$ and the automorphism group $\Aut(G)$ of $G$ are permutation groups on $G$. It is easy to see that $R(g)^\a=R(g^\a)$ for every $g\in G$ and $\a\in \Aut(G)$, implying $R(G)$ is normalized by $\Aut(G)$. Furthermore, $R(G)\cap \Aut(G)=1$, and hence we have $R(G)\Aut(G)=R(G)\rtimes\Aut(G)$. This group is called the {\em holomorph} of the group $G$, denoted by $\Hol(G)$, that is, $\Hol(G)=\langle R(G),\Aut(G)\rangle=R(G)\rtimes\Aut(G)$ (see \cite[Lemma 7.16]{Rotman}).

\medskip
Now we are ready to prove Theorem~\ref{mainth}.

\medskip

\noindent{\bf Proof of Theorem~\ref{mainth}:} First we prove the first half of Theorem~\ref{mainth}, that is, $C_n$ is a NDCI-group if and only if $8\nmid n$.

Let $8\mid n$. By Proposition \ref{pro-no-ci} and Lemma \ref{thm-no-ci}, there exists a normal non-CI Cayley digraph on $C_n$, that is, $C_n$ is not a NDCI-group. The necessity follows. To prove the sufficiency, assume that $8\nmid n$ and we only need to  prove that $C_n$ is a NDCI-group.

For convenience, write $n=p_1^{r_1}p_2^{r_2}\cdots p_m^{r_m}2^s$ as the prime factorization of $n$, where $s\leq 2$ and $p_1, p_2, \cdots, p_m$ are distinct odd prime factors of $n$ such that $p_1>p_2>\cdots>p_m$. By Proposition \ref{AdNe-DCI-group}, $C_n$ is a DCI-group for $n=1,2,4$, and so a NDCI-group. Thus, we may assume that $m\geq 1$. In particular, $r_i\geq 1$ for each $1\leq i \leq m$.

Let $\Gamma=\Cay(C_n,S)$ be a normal Cayley digraph and let $A=\Aut(\Gamma)$. To prove that $C_n$ is a NDCI-group, it suffices to show that $\Gamma$ is a CI-digraph.

Let $C_n=C_{p_1^{r_1}}\times C_{p_2^{r_2}}\times\cdots C_{p_m^{r_m}}\times C_{2^s}$. By Proposition~\ref{AGCG}, $$\Aut(C_n)=\Aut(C_{p_1^{r_1}})\times \Aut(C_{p_2^{r_2}})\times\cdots\times\Aut(C_{p_m^{r_m}})\times \Aut(C_{2^s}).$$ Recall that $\Hol(C_n)=R(C_n)\Aut(C_n)$, and $R(x)^\a=R(x^\a)$ for all $\a\in \Aut(C_n)$ and $x\in C_n$.  It is easy to see that for all $\a_i\in \Aut(C_{p_i^{r_i}})$,
$$(z_1\cdots z_{i-1} z_iz_{i+1}\cdots z_n)^{\a_i}=z_1\cdots z_{i-1}z_i^{\a_i}z_{i+1}\cdots z_m, \mbox{ where } z_j\in C_{p_j^{r^j}}.$$ It follows that for all $i\not=j$,  $$[R(C_{p_j^{r_j}}),\Aut(C_{p_i^{r_i}})]=[R(C_{2^s}),\Aut(C_{p_i^{r_i}})]=1,$$ that is, $\Aut(C_{p_i^{r_i}})$ commutes with $R(C_{p_j^{r_j}})$ and $R(C_{2^s})$ elementwise. Since $\Gamma=\Cay(C_n,S)$ is normal, Proposition~\ref{N_AUT} implies that $A=R(C_n)\Aut(C_n,S)\leq \Hol(C_n)$ and $A_1=\Aut(C_n,S)$.

Let $G$ be a regular subgroup of $A$ with $G\cong C_n$. To prove that $\Gamma$ is a CI-digraph, by Proposition~\ref{CI-graph-prop} it suffices to show that $G=R(C_n)$.

Note that $G$ has a unique Sylow $p_i$-subgroup $G_{p_i}$ with $|G_{p_i}|=p_i^{r_i}$, and $G\leq A\leq \Hol(C_n)=R(C_n)\rtimes (\Aut(C_{p_1^{r_1}})\times \Aut(C_{p_2^{r_2}})\times\cdots\times\Aut(C_{p_m^{r_m}})\times \Aut(C_{2^s}))$.

\medskip
\noindent {\bf Claim 1:} Assume that $g=R(a)\a_1\a_2\cdots\a_m\a_{m+1} \in G$, where $a\in C_n$, $\a_{m+1}\in \Aut(C_{2^s})$ and $\a_i\in \Aut(C_{p_i^{r_i}})$ for each $1\leq i\leq m$, and assume that $R(C_{p_k^{r_k}})\leq G$ for some $1\leq k\leq m$. Then $\a_k=1$.

Since $R(C_{p_k^{r_k}})\leq G$ and $G\cong C_n$, we have $[g,R(C_{p_k^{r_k}})]=1$. Since $[\Aut(C_{p_i^{r_i}}), R(C_{p_k^{r_k}})]=1$ for every $i\not=k$ and $[\Aut(C_{2^s}), R(C_{p_k^{r_k}})]=1$, we have $[\a_i,R(C_{p_k^{r_k}})]=1$ for every $i\not=k$, and since $[R(a),R(C_{p_k^{r_k}})]=1$ and $g=R(a)\a_1\a_2\cdots\a_m\a_{m+1}$, we have $[\a_k, R(C_{p_k^{r_k}})]=1$. This implies that $\a_k=1$ as $\a_k\in \Aut(C_{p_k^{r_k}})$.

\medskip

\noindent {\bf Claim 2:} $R(C_{p_k^{r_k}})\leq G$ for each $1\leq k\leq m$.

Recall that $\Hol(C_n)_{p_i}$ is a Sylow $p_i$-subgroup of $\Hol(C_n)$ and $R(C_n)\unlhd \Hol(C_n)=R(C_n)\Aut(C_n)$. Then $\Hol(C_n)_{p_i}\leq R(C_n)(\Aut(C_n))_{p_i}$ for every $1\leq i\leq m$, and since $\Aut(C_n)$ is abelian, $(\Aut(C_n))_{p_i}$ is the unique Sylow $p_i$-subgroup of $\Aut(C_n)$. Since $\Aut(C_{p_i^{r_i}})\cong C_{(p_i-1)p_i^{r_i-1}}$, $\Aut(C_{p_i^{r_i}})$ has a unique subgroup of order $p_i-1$, denoted by $\Aut(C_{p_i^{r_i}})_{p_i-1}$, and a unique Sylow $p_i$-subgroup, that is, $\Aut(C_{p_i^{r_i}})_{p_i}$. It follows that $\Aut(C_{p_i^{r_i}})=\Aut(C_{p_i^{r_i}})_{p_i-1}\times \Aut(C_{p_i^{r_i}})_{p_i}$. Since $p_1>p_2>\cdots>p_m$, we have $$\Aut(C_n)_{p_i}\leq \Aut(C_{p_1^{r_1}})_{p_1-1}\times \Aut(C_{p_2^{r_2}})_{p_2-1}\times \cdots\times \Aut(C_{p_{i-1}^{r_{i-1}}})_{p_{i-1}-1}\times \Aut(C_{p_i^{r_i}})_{p_i},$$ and in particular, $\Aut(C_n)_{p_1}= \Aut(C_{p_1^{r_1}})_{p_1}$.

We process the proof by induction on $k$. For simplicity, we may assume that $0\leq k\leq m$ and for $k=0$, let $R(C_{p_k^{r_k}})=1$. Then the claim is true for $k=0$.

Let $k\geq 1$. By inductive hypothesis, we assume that $R(C_{p_i^{r_i}})\leq G$ for every $0\leq i<k$. We suppose that  $R(C_{p_k^{r_k}})\nleq G$ and will obtain a contradiction. Since $|G_{p_k}|=p_k^{r_k}$, there exists $h\in G_{p_k}$ but $h\not\in R(C_{p_k^{r_k}})$. Since $h\in \Hol(C_n)_{p_k}\leq R(C_n)(\Aut(C_n))_{p_k}$ and $\Aut(C_n)_{p_k}\leq \Aut(C_{p_1^{r_1}})_{p_1-1}\times \Aut(C_{p_2^{r_2}})_{p_2-1}\times\cdots\times \Aut(C_{p_{k-1}^{r_{k-1}}})_{p_{k-1}-1}\times \Aut(C_{p_k^{r_k}})_{p_k}$, we have $h=R(h_1)\b$ for some $h_1\in C_n$ and $\b\in (\Aut(C_n))_{p_k}$, where $\b=\b_1\b_2\cdots \b_k$ for some $\b_j\in  \Aut(C_{p_j^{r_j}})_{p_j-1}$ with $1\leq j<k$ and $\b_k\in \Aut(C_{p_k^{r_k}})_{p_k}$. Note that $h$ is a $p_k$-element and $R(C_{p_k^{r_k}})$ is the unique Sylow $p_k$-subgroup of $R(C_n)$. Then $h\not\in R(C_{p_k^{r_k}})$ implies that $h\not\in R(C_n)$, forcing $\b\not=1$. Since $R(C_{p_i^{r_i}})\leq G$ for every $1\leq i<k$, by Claim~1 we have $\b_i=1$. It follows $\b=\b_k\not=1$, and $\b_k=R(h_1)^{-1}h\in R(C_n)G\leq A$. Furthermore, $\b_k\in A_1=\Aut(C_n,S)$, and hence $\Aut(C_n,S)$ contains an element of order $p_k$ in $\Aut(C_{p_k^{r_k}})$. By Lemma~\ref{non-normal}, $\Gamma$ is non-normal, a contradiction. Thus, $R(C_{p_k^{r_k}})\leq G$. By induction, $R(C_{p_k^{r_k}})\leq G$ for each $0\leq k\leq m$, as claimed.

\medskip
To prove $G=R(C_n)=R(C_{p_1^{r_1}})\times R(C_{p_2^{r_2}})\times\cdots\times R(C_{p_m^{r_m}})\times R(C_{2^s})$ with $s\leq 2$, by Claim~2 we are only left to show  $R(C_{2^s})\leq G$, and it suffices to show that $G_2\leq R(C_{2^s})$, because $|G_2|=2^s$.

Let $x\in G_2$. Since $$G\leq A\leq \Hol(C_n)=R(C_n)\rtimes (\Aut(C_{p_1^{r_1}})\times \Aut(C_{p_2^{r_2}})\times\cdots\times\Aut(C_{p_m^{r_m}})\times \Aut(C_{2^s})),$$ we have $x=R(x_1)\g_1\g_2\cdots\g_m\a$, where $x_1\in C_n$, $\a\in \Aut(C_{2^s})$ and $\g_i\in  \Aut(C_{p_i^{r_i}})_{p_i-1}$ for $1\leq i\leq m$. By Claims~2 and 1, $\g_i=1$ for each $1\leq i\leq m$, and hence $x=R(x_1)\a\in R(C_n)\Aut(C_{2^s})$. Since  $[R(C_{p_i^{r_i}}),\Aut(C_{2^s})]=1$, $R(C_n)\Aut(C_{2^s})$ has a unique Sylow $2$-subgroup, that is, $R(C_{2^s})\Aut(C_{2^s})$. It follows that  $x\in R(C_{2^s})\Aut(C_{2^s})$, and hence $G_2\leq R(C_{2^s})\Aut(C_{2^s})$.

If $\Aut(C_{2^s})=1$, then $G_2\leq R(C_{2^s})$, as required. We may assume that $\Aut(C_{2^s})\not=1$, and hence $s=2$ as $s\leq 2$. It follows that $\Aut(C_{2^s})=\Aut(C_4)\cong C_2$. In particular, $R(C_{2^s})\Aut(C_{2^s})$ is a dihedral group of order $8$ and hence has a unique cyclic subgroup of order $4$, that is, $R(C_{2^s})$. Since $G\cong C_n$, $G_2$ is a cyclic group of order $4$, and since $G_2\leq R(C_{2^s})\Aut(C_{2^s})$, we obtain that $G_2=R(C_{2^s})$, as required. This completes the proof of the first half of the theorem.

Now we are ready to prove the second half of the theorem, that is,  $C_n$ is a $\rm{NCI}$-group if and only if either $n=8$ or $8\nmid n$. If $n=8$ then by Proposition~\ref{AdNe-DCI-group}, $C_n$ is a CI-group and hence a NCI-group; if $8\nmid n$ then by above proof of the first half of Theorem~\ref{mainth}, $C_n$ is a NDCI-group and hence a NCI-group. This proves the sufficiency, and the necessity follows from Lemma~\ref{thm-no-ci}.   \qed

\section{Ideas from S-ring theory}

One of the approaches to study CI-property is based on the method of S-rings. This approach was suggested in~\cite{HM} and realized in~\cite{Feng,KM,Kov,MS}. A Cayley digraph $\Gamma=\Cay(G,S)$ of a group $G$ is CI if and only if the smallest S-ring $\mathcal{A}$ over $G$ such that $\underline{S}=\sum_{s\in S} s\in \mathcal{A}$ is CI (see~\cite[Definition~3]{HM} for the definition of a CI-S-ring). The automorphism group of any S-ring over $G$ contains the group $R(G)$ and an S-ring is called {\em normal} if $R(G)$ is normal in its automorphism group. Due to~\cite[Theorem~2.6.4]{CP}, we have $\Aut(\Gamma)=\Aut(\mathcal{A})$. So $\Gamma$ is normal if and only if $\mathcal{A}$ is normal. Normal S-rings over cyclic groups were studied in~\cite{EP}. In fact, using~\cite[Corollary~6.5, Theorem~6.6, Lemma~7.1]{EP}, it is possible to derive that if $8\nmid n$ then every normal S-ring and hence every normal Cayley digraph over a cyclic group of order~$n$ is CI.
\vspace{5mm}

\noindent {\bf Acknowledgements:} The first, second, and fourth authors were supported by the National Natural Science Foundation of China (11731002) and the 111 Project of China (B16002). The third author was supported by the Mathematical Center in Akademgorodok under the agreement No. 075-15-2019-1613 with the Ministry of Science and Higher Education of the Russian Federation.

\end{document}